\documentclass{ifacconf}

\usepackage[utf8]{inputenc}
\usepackage{amsmath}
\usepackage{enumerate}
\usepackage{amsfonts}
\usepackage{amssymb}
\usepackage{graphicx}

\usepackage{subcaption}
\usepackage{color}
\usepackage{enumitem}
\usepackage{amsbsy}
\usepackage{bm}
\newcommand{\bmu}{\bm{u}}
\newcommand{\bmx}{\bm{x}}

\newcommand{\bmw}{\bm{w}}

\newcommand{\bmz}{\bm{z}}

\DeclareMathOperator{\is}{i/s}

\DeclareMathOperator{\ids}{i/d/s}
\DeclareMathOperator{\im}{im}

\DeclareMathOperator{\rank}{rank}
\DeclareMathOperator{\trace}{tr}

\let\leq\leqslant
\let\geq\geqslant

\newcommand{\bmat}{\begin{matrix}}
\newcommand{\emat}{\end{matrix}}
\newcommand{\bbm}{\begin{bmatrix}}
\newcommand{\ebm}{\end{bmatrix}}
\newcommand{\bbma}{\begin{bmatrix*}[r]}
\newcommand{\ebma}{\end{bmatrix*}}
\newcommand{\bpm}{\begin{pmatrix}}
\newcommand{\epm}{\end{pmatrix}}
\newcommand{\bvm}{\begin{vmatrix}}
\newcommand{\evm}{\end{vmatrix}}
\newcommand{\bse}{\begin{subequations}}
\newcommand{\ese}{\end{subequations}}
\newcommand{\beq}{\begin{equation}}
\newcommand{\eeq}{\end{equation}}
\newcommand{\ben}{\renewcommand{\labelenumi}{\arabic{enumi}.}
\renewcommand{\theenumi}{\arabic{enumi}}\begin{enumerate}}

\newcommand{\een}{\end{enumerate}}

\newcommand{\beni}{\renewcommand{\labelenumi}{\roman{enumi}.}
\renewcommand{\theenumi}{\roman{enumi}}\begin{enumerate}}

\newcommand{\eeni}{\end{enumerate}}

\newcommand{\bena}{\renewcommand{\labelenumi}{\alph{enumi}.}
\renewcommand{\theenumi}{\alph{enumi}}\begin{enumerate}}

\newcommand{\eena}{\end{enumerate}}

\newcommand{\bit}{\begin{itemize}}
\newcommand{\eit}{\end{itemize}}
\newcommand{\bthe}{\begin{theorem}}
\newcommand{\ethe}{\end{theorem}}
\newcommand{\blem}{\begin{lemma}}
\newcommand{\elem}{\end{lemma}}
\newcommand{\bprop}{\begin{proposition}}
\newcommand{\eprop}{\end{proposition}}
\newcommand{\bex}{\begin{example}}
\newcommand{\eex}{\end{example}}
\newcommand{\bas}{\begin{assumption}}
\newcommand{\eas}{\end{assumption}}
\newcommand{\bre}{\begin{remark}}
\newcommand{\ere}{\end{remark}}
\newcommand{\bcor}{\begin{corollary}}
\newcommand{\ecor}{\end{corollary}}
\newcommand{\bdfn}{\begin{definition}}
\newcommand{\edfn}{\end{definition}}
\newcommand{\bcon}{\begin{conjecture}}
\newcommand{\econ}{\end{conjecture}}

\usepackage{graphicx}      % include this line if your document contains figures
\usepackage{natbib}        % required for bibliography
\begin{document}
\begin{frontmatter}

\title{Data-driven parameterizations of suboptimal LQR and $\mathcal{H}_2$ controllers}

\thanks[Henk]{The first author acknowledges financial support by the RAIN lab at University of Washington and the Centre for Data Science and Systems Complexity at University of Groningen.}

\author[Henk]{Henk J. van Waarde} 
\author[Mehran]{Mehran Mesbahi} 

\address[Henk]{Bernoulli Institute for Mathematics, Computer Science and Artificial Intelligence, and Engineering and Technology Institute Groningen, University of Groningen, the Netherlands 
   (h.j.van.waarde@rug.nl)}
\address[Mehran]{William E. Boeing Department of Aeronautics and Astronautics, University of Washington, Seattle, WA 98195 USA (mesbahi@uw.edu)}

\begin{abstract}                
In this paper we design suboptimal control laws for an unknown linear system on the basis of measured data. We focus on the suboptimal linear quadratic regulator problem and the suboptimal $\mathcal{H}_2$ control problem. For both problems, we establish conditions under which a given data set contains sufficient information for controller design. We follow up by providing a data-driven parameterization of all suboptimal controllers. We will illustrate our results by numerical simulations, which will reveal an interesting trade-off between the number of collected data samples and the achieved controller performance. 
\end{abstract}


\end{frontmatter}
%===============================================================================

\section{Introduction}
In the field of systems and control, the majority of control techniques is \emph{model-based}, meaning that these methods require knowledge of a plant model, for example in the form of a transfer function or state-space system. Such system models are rarely known a priori and typically have to be identified using measured data. The aim of \emph{data-driven} control is to bypass this system identification step, and to design control laws for dynamical systems directly on the basis of data. Contributions to data-driven control can roughly be divided in on- and offline techniques. 

Methods in the former class are iterative and make use of multiple online experiments. Examples include direct adaptive control (\cite{Astrom1989}), iterative feedback tuning (\cite{Hjalmarsson1998}) and methods based on reinforcement learning (\cite{Bradtke1993,Alemzadeh2019}). Offline techniques construct controllers on the basis of data (typically a single system trajectory) that is collected offline. \cite{Skelton1994} consider optimal control using a batch-form solution to the Riccati equation. Virtual reference feedback tuning was introduced by \cite{Campi2002}. Moreover, \cite{Campestrini2017} cast the problem of designing model reference controllers in the prediction error framework. \cite{Baggio2019} design minimum energy controls using data. The fundamental lemma by \cite{Willems2005} has also been leveraged for data-driven control in a behavioral setting (\cite{Markovsky2008}), and in the context of state-space systems to design model predictive controllers (\cite{Coulson2019}), stabilizing and optimal controllers (\cite{DePersis2019}) and robust controllers (\cite{Berberich2019c}).

An important persisting problem is to understand the relative merits of data-driven control and combined system identification and model-based control, see e.g. (\cite{Tu2018}). A recent paper sheds some light on this issue by studying data-driven control from the perspective of \emph{data informativity}. In particular, \cite{vanWaarde2019} provide conditions under which given data contain enough information for control design. For control problems such as stabilization, these conditions do not require that the underlying system can be uniquely identified. As such, one can generally stabilize an unknown system without learning its dynamics exactly. For the linear quadratic regulator problem, however, it was shown that the data essentially need to be rich enough for system identification. 

Inspired by the above results, it is our goal to study data-driven \emph{suboptimal} control problems. Intuitively, we expect that the data requirements for such suboptimal problems are \emph{weaker} than those for their optimal counterparts. We will focus on data-driven versions of the suboptimal linear quadratic regulator (LQR) problem and the $\mathcal{H}_2$ suboptimal control problem. Both of these problems involve the data-guided design of controllers that stabilize the unknown system and render the (LQR or $\mathcal{H}_2$) cost smaller than a given tolerance. 

Our main results are the following. First, for both suboptimal problems, we establish necessary and sufficient conditions under which the data are informative for control design. These conditions do not require that the underlying system can be identified uniquely. Secondly, for both problems we give a parameterization of all suboptimal controllers in terms of data-driven linear matrix inequalities.

\emph{Outline:} In \S \ref{sectionsuboptimal} we provide some preliminaries. In \S \ref{sectionproblem} we state the problem. Next, \S \ref{sectionLQR} and \S \ref{sectionH2} contain our main results. An illustrative example is given in \S \ref{sectionexample}. Finally, \S \ref{sectionconclusions} contains our conclusions. 

\newpage

\section{Suboptimal control problems}
\label{sectionsuboptimal}
The purpose of this section is to review two (model-based) suboptimal control problems whose data-driven versions will be the main topic of this paper. 

\subsection{The suboptimal LQR problem}
Consider the linear system 
\begin{equation}
\label{sysAB}
    \bmx(t+1) = A\bmx(t) + B\bmu(t),
\end{equation}
where $\bmx \in \mathbb{R}^n$ is the state, $\bmu \in \mathbb{R}^m$ is the input and $A$ and $B$ are real matrices of appropriate dimensions. We will occasionally use the shorthand notation $(A,B)$ to refer to system \eqref{sysAB}. Associated with \eqref{sysAB}, we consider the infinite-horizon cost functional
\begin{equation}
\label{cost}
J(x_0,u)=\sum_{t=0}^\infty  x^\top(t) Q x(t) + u^\top(t) R u(t),
\end{equation}
where $x_0$ is the initial state and $Q = Q^\top \geq 0$ and $R = R^\top > 0$ are real matrices. Whenever the input function $u$ results from a state feedback law $\bmu = K\bmx$, we will write $J(x_0,K)$ instead of $J(x_0,u)$. The suboptimal linear quadratic regulator problem can be formulated as follows. Given an initial condition $x_0 \in \mathbb{R}^n$ and tolerance $\gamma > 0$, find (if it exists) a feedback law $\bmu = K\bmx$ such that $A+BK$ is stable\footnote{Here we refer to the notion of \emph{Schur stability}, i.e., a matrix is said to be stable if all its eigenvalues are contained in the open unit disk.}, and the cost satisfies $J(x_0,K) < \gamma$. Such a $K$ is called a \emph{suboptimal feedback gain} for the system $(A,B)$. The following proposition gives necessary and sufficient conditions under which a given matrix $K$ is a suboptimal feedback gain.
\begin{prop}
\label{propsuboptimalLQR}
Let $x_0 \in \mathbb{R}^n$ and $\gamma > 0$. The matrix $K$ is a suboptimal feedback gain if and only if there exists a matrix $P = P^\top > 0$ such that 
\begin{align}
    (A+BK)^\top P (A+BK) - P + Q + K^\top R K &< 0 \label{eq1} \\
    x_0^\top P x_0 &< \gamma. \label{eq2}
\end{align}
\end{prop}

\subsection{The $\mathcal{H}_2$ suboptimal control problem}

Consider the system
\begin{subequations}
\label{ABECD}
\begin{align}
    \label{ABE}
    \bmx(t+1) &= A\bmx(t) + B\bmu(t) + E\bmw(t) \\
    \bmz(t) &= C\bmx(t) + D\bmu(t), \label{CD}
    \end{align} 
\end{subequations}
where $\bmx \in \mathbb{R}^n$ denotes the state, $\bmu \in \mathbb{R}^m$ is the control input, $\bmw \in \mathbb{R}^d$ is a disturbance input and $\bmz \in \mathbb{R}^p$ is the performance output. The real matrices $A,B,C,D$ and $E$ are of appropriate dimensions. The feedback law $\bmu = K \bmx$ yields the closed-loop system 
\begin{subequations}
\label{closedloop}
\begin{align}
\label{closedloop1}
    \bmx(t+1) &= (A+BK)\bmx(t) + E\bmw(t) \\
    \label{closedloop2}
    \bmz(t) &= (C+DK)\bmx(t).
\end{align}
\end{subequations}
Associated with \eqref{closedloop}, we consider the $\mathcal{H}_2$ cost functional
$$
J_{\mathcal{H}_2}(K) := \sum_{t=0}^\infty \trace\left( T_K^\top(t) T_K(t) \right),
$$
where $T_K(t) := (C+DK)(A+BK)^t E$ is the closed-loop impulse response from $\bmw$ to $\bmz$ and $\trace$ denotes trace. The cost $J_{\mathcal{H}_2}(K)$ equals the squared $\mathcal{H}_2$ norm of the transfer function from $\bmw$ to $\bmz$ of \eqref{closedloop}. It is well-known that the $\mathcal{H}_2$ cost of a given stabilizing $K$ can be computed using the observability Gramian. Indeed for a stabilizing $K$, the unique solution $P$ to the Lyapunov equation
\begin{equation} 
\label{Lyapunoveq}
(A+BK)^\top P (A+BK) - P + (C+DK)^\top (C+DK) = 0
\end{equation} 
is related to the $\mathcal{H}_2$ cost by $\trace(E^\top P E) = J_{\mathcal{H}_2}(K)$. For a given $\gamma > 0$, the $\mathcal{H}_2$ suboptimal control problem amounts to finding a gain $K$ (if it exists) such that $A+BK$ is stable and $J_{\mathcal{H}_2}(K) < \gamma$. Such a $K$ is called an \emph{$\mathcal{H}_2$ suboptimal feedback gain}. Similar to Proposition \ref{propsuboptimalLQR} the following proposition gives conditions under which a given $K$ is an $\mathcal{H}_2$ suboptimal feedback gain. 
\begin{prop}
\label{propsuboptimalH2}
Let $\gamma > 0$. The matrix $K$ is an $\mathcal{H}_2$ suboptimal feedback gain if and only if there exists a matrix $P = P^\top > 0$ such that 
\begin{align*}
    (A+BK)^\top P (A+BK) - P + (C+DK)^\top (C+DK) &< 0 \\
    \trace(E^\top P E) &< \gamma.
\end{align*}
\end{prop}
Clearly, the LQR suboptimal control problem can be viewed as a special case of the $\mathcal{H}_2$ suboptimal control problem. Indeed, the $\mathcal{H}_2$ problem boils down to the LQR problem if $E = x_0$, $C^\top C = Q$, $D^\top D = R$ and $C^\top D = 0$. However, as we will see in the next section, the data-driven versions of these problems are different in the way that data is collected.

\section{Problem formulation}
\label{sectionproblem}
In this section we formulate our problems. We will start by introducing the data-driven suboptimal LQR problem. Consider the linear system
\begin{equation}
\label{system}
\bmx(t+1) = A_s \bmx(t) + B_s \bmu(t),
\end{equation}
where $\bmx \in \mathbb{R}^n$ denotes the state, $\bmu \in \mathbb{R}^m$ is the input and $A_s$ and $B_s$ are real matrices of appropriate dimensions. We refer to \eqref{system} as the `true' system. Suppose that the system matrices $A_s$ and $B_s$ of the true system are unknown, but we have access to a finite set of data\footnote{We assume a single trajectory is measured. Our results are also applicable in case multiple (short) trajectories are measured, which can be beneficial if $A_s$ is unstable \citep{vanWaarde2020}.}
    \begin{align*}
    	U_- &:= \bbm u(0) & u(1) & \cdots & u(T-1)\ebm \\
    	X_{\phantom{-}} &:= \bbm x(0) & x(1) & \cdots & x(T)\ebm,
    \end{align*}
generated by system \eqref{system}. By partitioning the state data as 
    \begin{align*}
		X_- &:= \bbm x(0) & x(1) & \cdots & x(T-1) \ebm \\
		X_+ &:= \bbm x(1) & x(2) & \cdots & x(T) \ebm,
	\end{align*}
we can relate the data and $(A_s,B_s)$ via
    \begin{equation*}
        X_+ = \begin{bmatrix} A_s & B_s \end{bmatrix} \begin{bmatrix} X_- \\ U_- \end{bmatrix}.
    \end{equation*}
The set of all systems that explain the input/state data $(U_-,X)$ is given by    
    \begin{equation*}
	\Sigma_{\is} := \left\{ (A,B) \mid X_+ = \bbm A&B \ebm
	\begin{bmatrix}
	X_-\\U_-
	\end{bmatrix} \right\}.
	\end{equation*} 
Associated with system \eqref{system} we consider the cost functional \eqref{cost}, where the matrices $Q = Q^\top \geq 0$ and $R = R^\top > 0$ and the initial condition\footnote{We emphasize that the initial condition $x_0$ is not necessarily the same as the first measured state sample $x(0)$.} $x_0$ are assumed to be given. We want to design a suboptimal feedback gain for the unknown $(A_s,B_s)$ on the basis of the data. Given $(U_-,X)$, it is impossible to distinguish between the systems in $\Sigma_{\is}$, and therefore we can only guarantee that $K$ is a suboptimal gain for $(A_s,B_s)$ if it is a suboptimal gain \emph{for all} systems in $\Sigma_{\is}$. With this in mind, we introduce the following notion of data informativity.
\begin{defn}
Let $x_0 \in \mathbb{R}^n$ and $\gamma > 0$. The data $(U_-,X)$ are \emph{informative for suboptimal linear quadratic regulation} if there exists a matrix $K$ that is a suboptimal feedback gain for all $(A,B) \in \Sigma_{\is}$.
\end{defn}
We want to find conditions under which the data are informative for suboptimal LQR, and we want to obtain suboptimal controllers from data. These problems are stated more formally as follows.

\begin{prob}
\label{problemLQR}
Let $x_0 \in \mathbb{R}^n$ and $\gamma > 0$. Provide necessary and sufficient conditions under which the data $(U_-,X)$ are informative for suboptimal linear quadratic regulation. Moreover, for data $(U_-,X)$ that are informative, find a feedback gain $K$ that is suboptimal for all $(A,B) \in \Sigma_{\is}$.
\end{prob}

Subsequently, we turn our attention to the $\mathcal{H}_2$ suboptimal control problem. For this, consider the system
\begin{align}
    \bmx(t+1) &= A_s\bmx(t) + B_s\bmu(t) + E_s\bmw(t) \\
    \bmz(t) &= C\bmx(t) + D\bmu(t),
\end{align}
where the system matrices $A_s$, $B_s$ and $E_s$ are unknown, but the matrices $C$ and $D$ defining the performance output are known. We collect the data $X$ and $U_-$ as before, as well as the corresponding measurements of the disturbance 
$$
W_- := \bbm w(0) & w(1) & \cdots & w(T-1)\ebm. 
$$ 
The assumption that $W_-$ is available is reasonable in applications such as aircraft control, where gust disturbances can be measured via on-board LIDAR measurement systems, see e.g., \cite{Soreide1996}.
In this setup, all triples of system matrices $(A,B,E)$ that explain the data $(U_-,W_-,X)$ are given by 
\begin{equation*}
    \Sigma_{\ids} := \left\{ (A,B,E) \mid X_+ = \begin{bmatrix} A & B & E \end{bmatrix} \begin{bmatrix} X_- \\ U_- \\ W_- \end{bmatrix} \right\}.
\end{equation*}

We can now state the following notion of data informativity for $\mathcal{H}_2$ suboptimal control.

\begin{defn}
Let $\gamma > 0$. The data $(U_-,W_-,X)$ are \emph{informative for $\mathcal{H}_2$ suboptimal control} if there exists a $K$ that is an $\mathcal{H}_2$ suboptimal feedback gain for all $(A,B,E) \in \Sigma_{\ids}$.
\end{defn}

As before, we are interested in both data informativity conditions and a control design procedure. We formalize this in the following problem. 
\begin{prob}
\label{problemH2}
Let $\gamma > 0$. Provide necessary and sufficient conditions under which the data $(U_-,W_-,X)$ are informative for $\mathcal{H}_2$ suboptimal control. Moreover, for data $(U_-,W_-,X)$ that are informative, find a feedback gain $K$ that is $\mathcal{H}_2$ suboptimal for all $(A,B) \in \Sigma_{\is}$.
\end{prob}

\begin{rem}
We note that the data-driven $\mathcal{H}_2$ \emph{optimal} control problem was studied by \cite{DePersis2019} in the case that $E_s = I$ and $(U_-,X)$ data are collected in the absence of disturbances. Sufficient data conditions were given for this problem via the concept of persistency of excitation. Moreover, \cite{Berberich2019c} aim to design data-driven controllers that minimize a quadratic performance specification (with the $\mathcal{H}_{\infty}$ problem as a special case). The authors provide sufficient data conditions in the scenario that $E$ is known and $\bmw$ is unmeasured. 
\end{rem}

\section{Data-driven suboptimal LQR}
\label{sectionLQR}
In this section we report our solution to Problem \ref{problemLQR}. Before we start, we need some results from (\cite{vanWaarde2019}). We say that $(U_-,X)$ are \emph{informative for stabilization by state feedback} if there exists a $K$ such that $A+BK$ is stable for all $(A,B) \in \Sigma_{\is}$. The following result was proven in \cite[Thm. 16]{vanWaarde2019}. 
\begin{lem}
\label{leminfstab}
The data $(U_-,X)$ are informative for stabilization by state feedback if and only if there exists a right inverse $X_-^\dagger$ of $X_-$ such that $X_+ X_-^\dagger$ is stable. 
	
Moreover, $K$ is a stabilizing feedback for all systems in $\Sigma_{\is}$ if and only if $K = U_- X_-^\dagger$ for some $X_-^\dagger$ satisfying the above properties.
\end{lem}

Next, we characterize the informativity of data for suboptimal LQR in terms of data-driven matrix inequalities. 

\begin{thm}
\label{thmineq}
Let $x_0 \in \mathbb{R}^n$ and $\gamma > 0$. The data $(U_-,X)$ are informative for suboptimal linear quadratic regulation if and only if there exists a matrix $P = P^\top > 0$ and a right inverse $X_-^\dagger$ of $X_-$ such that 
\begin{flalign}
\label{ineq}
    (X_+X_-^\dagger)^\top P X_+X_-^\dagger - P + Q + (U_- X_-^\dagger)^\top R U_- X_-^\dagger &<0&& \\
    x_0^\top P x_0 &<\gamma.&& \label{ineqcost}
\end{flalign}    
Moreover, $K$ is a suboptimal feedback gain for all systems $(A,B) \in \Sigma_{\is}$ if and only if it is of the form $K = U_- X_-^\dagger$ for some right inverse $X_-^\dagger$ satisfying \eqref{ineq} and \eqref{ineqcost}.
\end{thm}

\begin{pf}
To prove the `if' parts of both statements, suppose that there exists a matrix $P = P^\top > 0$ and a right inverse $X_-^\dagger$ such that \eqref{ineq} and \eqref{ineqcost} are satisfied. Define the controller $K := U_- X_-^\dagger$. For any $(A,B) \in \Sigma_{\is}$ we have 
$
X_+ = AX_- + BU_-,
$
which implies that $X_+ X_-^\dagger = A+BK$. Substitution of the latter expression into \eqref{ineq} yields \begin{align*}
    (A+BK)^\top P (A+BK) - P + Q + K^\top R K < 0,
\end{align*}
which shows that there exists a $K$ and $P = P^\top > 0$ satisfying \eqref{eq1} and \eqref{eq2} for all $(A,B) \in \Sigma_{\is}$. By Proposition \ref{propsuboptimalLQR}, the data are informative for suboptimal LQR. 

To prove the `only if' parts of both statements, suppose that the data $(U_-,X)$ are informative for suboptimal linear quadratic regulation. This means that there exists a feedback gain $K$ and a matrix $P_{(A,B)} = P_{(A,B)}^\top > 0$ such that 
\begin{align*}
    (A+BK)^\top P_{(A,B)} (A+BK) - P_{(A,B)} + Q + K^\top R K &< 0  \\
    x_0^\top P_{(A,B)} x_0 &< \gamma 
\end{align*}
for all $(A,B) \in \Sigma_{\is}$. We emphasize that the matrix $P_{(A,B)}$ may depend on the particular system $(A,B)$, but the feedback gain $K$ is fixed by definition. Since $K$ is such that $A+BK$ is stable for all $(A,B) \in \Sigma_{\is}$, we obtain by Lemma \ref{leminfstab} that $K$ is of the form $K = U_- X_-^\dagger$ for some right inverse $X_-^\dagger$ of $X_-$. This yields $A+BK = X_+ X_-^\dagger$. The matrix $A+BK$ is therefore the same for all $(A,B) \in \Sigma_{\is}$. This implies the existence of a (common) $P = P^\top > 0$ such that \eqref{ineq} and \eqref{ineqcost} are satisfied. \hfill $\square$
\end{pf}

Note that the conditions of Theorem \ref{thmineq} are not ideal from computational point of view since \eqref{ineq} depends nonlinearly on $P$ and $X_-^\dagger$. Nonetheless, it is straightforward to reformulate these conditions in terms of linear matrix inequalities. This is described in the following corollary. 

\begin{cor}
\label{corsuboptimalLQR}
Let $Q = C^\top C$, $R = D^\top D$ and $C^\top D = 0$, and let $x_0 \in \mathbb{R}^n$ and $\gamma > 0$. The data $(U_-,X)$ are informative for suboptimal linear quadratic regulation if and only if there exist $Y = Y^\top \in \mathbb{R}^{n \times n}$ and $\Theta \in \mathbb{R}^{T \times n}$ such that 
\begin{align}
    \begin{bmatrix}
    Y & \Theta^\top X_+^\top & \Theta^\top Z_-^\top \\
    X_+ \Theta & Y & 0 \\
    Z_-\Theta & 0 & I
    \end{bmatrix} &> 0 && \label{1} \\
    \begin{bmatrix}
    \gamma & x_0^\top \\ x_0 & Y
    \end{bmatrix} &> 0 && \label{2} \\
    X_-\Theta &= Y. && \label{3}
\end{align}
Here $Z_- := CX_- +DU_-$. Moreover, $K$ is a suboptimal feedback gain for all $(A,B) \in \Sigma_{\is}$ if and only if $K = U_- \Theta Y^{-1}$ for some $Y$ and $\Theta$ satisfying \eqref{1}, \eqref{2} and \eqref{3}.
\end{cor}

Corollary \ref{corsuboptimalLQR} follows from Theorem \ref{thmineq} via a few well-known tricks, see e.g. \cite{Scherer1999}. First a congruence transformation $P^{-1}$ is applied to \eqref{ineq}, after which a Schur complement argument and change of variables $Y := P^{-1}$ and $\Theta := X_-^\dagger Y$ yields \eqref{1}, \eqref{2} and \eqref{3}.

%\begin{pf}
%To prove the `if' part of both statements, suppose that there exists a symmetric $Z$ and matrix %$\Theta$ satisfying \eqref{1}, \eqref{2} and \eqref{3}. This is equivalent to the matrix
%\begin{align*}
%    Z^{-1} - Z^{-1}\Theta^\top (X_+^\top Z^{-1} X_+ + X_-^\top Q X_- + U_-^\top R U_-) \Theta %Z^{-1} 
%\end{align*}
%being positive definite and 
%\begin{align*} 
%    \gamma - x_0^\top Z^{-1} x_0 &> 0 \\
%    X_- \Theta Z^{-1} &= I,
%\end{align*}
%where we have applied the congruence transformation $Z^{-1}$ and used a Schur complement argument. Note that $\Theta Z^{-1}$ is a right inverse of $X_-$. Therefore, we see that \eqref{ineq} and \eqref{ineqcost} are satisfied for $P := Z^{-1}$ and $X_-^\dagger := \Theta Z^{-1}$. By Theorem \ref{thmineq}, the data are informative for suboptimal linear quadratic regulation. By the same lemma, the controller $K := U_-\Theta (X_- \Theta)^{-1}$ is a suboptimal feedback gain for all $(A,B) \in \Sigma_{\is}$.

%To prove the `only if' part of both statements, suppose that the data are informative for suboptimal linear quadratic regulation, and that $K$ is a suboptimal feedback gain. Then there exist matrices $P = P^\top > 0$ and $X_-^\dagger$ such that \eqref{ineq} and \eqref{ineqcost} are satisfied, and $K = U_- X_-^\dagger$. Using a Schur complement argument, we again see that $Z = P^{-1}$ and $\Theta = X_-^\dagger P^{-1}$ satisfy \eqref{1}, \eqref{2} and \eqref{3}, and $K = U_- \Theta (X_- \Theta)^{-1}$ for some $\Theta$ satisfying \eqref{1}, \eqref{2} and \eqref{3}. \hfill $\square$
%\end{pf}

\begin{rem} 
It is noteworthy that the conditions of Theorem \ref{thmineq} and Corollary \ref{corsuboptimalLQR} do not require that the data $(U_-,X)$ contain enough information to uniquely identify the system matrices $(A_s,B_s)$. Quite naturally, the conditions do become more difficult to satisfy for decreasing values of $\gamma$. Clearly, Theorem \ref{thmineq} and Corollary \ref{corsuboptimalLQR} require the matrix $X_-$ to have full row rank. This means that at least $T \geq n$ samples are needed to obtain a suboptimal controller from data. In comparison, note that to uniquely identify $A_s$ and $B_s$, it is necessary that the rank condition
$$
\rank 
\begin{bmatrix} X_- \\ U_- \end{bmatrix} = n+m
$$
is satisfied, which is only possible if $T \geq n + m$. In \S \ref{sectionexample} we will illustrate Corollary \ref{corsuboptimalLQR} in detail by numerical examples.
\end{rem}

\section{Data-driven $\mathcal{H}_2$ suboptimal control}
\label{sectionH2}

In this section we study the data-driven $\mathcal{H}_2$ suboptimal control problem as formulated in Problem \ref{problemH2}. As a first step, we extend Lemma \ref{leminfstab} to systems with disturbances. We say the data $(U_-,W_-,X)$ are \emph{informative for stabilization by state feedback} if there exists $K$ such that $A+BK$ is stable for all $(A,B,E) \in \Sigma_{\ids}$.

\begin{lem}
\label{lemstabE}
The data $(U_-,W_-,X)$ are informative for stabilization by state feedback if and only if there exists a right inverse $X_-^\dagger$ of $X_-$ with the properties that $X_+ X_-^\dagger$ is stable and $W_- X_-^\dagger = 0$. 

Moreover, $K$ is a stabilizing controller for all systems in $\Sigma_{\ids}$ if and only if $K = U_- X_-^\dagger$, where $X_-^\dagger$ satisfies the above properties. 
\end{lem}

\begin{pf}
The proof follows a similar line as that of \cite[Thm. 16]{vanWaarde2019}. To prove the `if' part of both statements, suppose that there exists a right inverse $X_-^\dagger$ such that $X_+ X_-^\dagger$ is stable and $W_- X_-^\dagger = 0$. Define $K := U_- X_-^\dagger$. Then $X_+ X_-^\dagger = A+BK$ for all $(A,B,E) \in \Sigma_{\ids}$. Hence $A+BK$ is stable for all $(A,B,E) \in \Sigma_{\ids}$ and $K = U_- X_-^\dagger$ is stabilizing.

To prove the `only if' parts, suppose that the data are informative for stabilization by state feedback. Let $K$ be stabilizing for all systems in $\Sigma_{\ids}$. Define the subspace
\begin{equation*}
    \Sigma_{\ids}^0 := \left\{ (A_0,B_0,E_0) \mid 0 = \begin{bmatrix} A_0 & B_0 & E_0 \end{bmatrix} \begin{bmatrix} X_- \\ U_- \\ W_- \end{bmatrix} \right\}.
\end{equation*}
The matrix $A+BK + \alpha (A_0+B_0K)$ is stable for all $\alpha \in \mathbb{R}$ and all $(A_0,B_0,E_0) \in \Sigma_{\ids}^0$. Thus we have
\begin{equation*}
    \rho\left( \frac{1}{\alpha} (A+BK) + A_0+B_0K \right) \leq \frac{1}{\alpha} \quad \forall \: \alpha \geq 1,
\end{equation*}
where $\rho(\cdot)$ denotes spectral radius. We take the limit as $\alpha \to \infty$, and conclude by continuity of the spectral radius that $A_0 + B_0 K$ is nilpotent for all $(A_0,B_0,E_0) \in \Sigma_{\ids}^0$. Note that $(A_0,B_0,E_0) \in \Sigma_{\ids}^0$ implies that 
$$
\left( (A_0 + B_0 K)^\top A_0, (A_0 + B_0 K)^\top B_0, (A_0 + B_0 K)^\top E_0 \right)
$$
is also a member of $\Sigma_{\ids}^0$. This implies that the matrix $(A_0 + B_0 K)^\top(A_0 + B_0 K)$ is nilpotent for all $(A_0,B_0,E_0)$. The only symmetric nilpotent matrix is zero, thus $A_0 + B_0 K = 0$ for all $(A_0,B_0,E_0) \in \Sigma_{\ids}^0$. We conclude that 
$$
\ker \begin{bmatrix} X_-^\top & U_-^\top & W_-^\top \end{bmatrix} \subseteq \ker \begin{bmatrix} I & K^\top & 0 \end{bmatrix},
$$
equivalently, 
$$
\im \begin{bmatrix}
I & K^\top & 0
\end{bmatrix}^\top \subseteq \im \begin{bmatrix} X_-^\top & U_-^\top & W_-^\top \end{bmatrix}^\top.
$$
This means that there exists a right inverse $X_-^\dagger$ of $X_-$ such that $K = U_- X_-^\dagger$ and $W_- X_-^\dagger = 0$. Clearly, $X_+ X_-^\dagger = A+BK$ for all $(A,B,E) \in \Sigma_{\ids}$, hence $X_+ X_-^\dagger$ is stable. \hfill $\square$
\end{pf}
The following theorem provides necessary and sufficient conditions for data informativity for the $\mathcal{H}_2$ problem. It also characterizes all suboptimal controllers in terms of the data. Recall that $Z_-$ was defined as $Z_- = CX_- + DU_-$.
\begin{thm}
\label{thmH2}
Let $\gamma > 0$. The data $(U_-,W_-,X)$ are informative for $\mathcal{H}_2$ suboptimal control if and only if at least one of the following two conditions is satisfied:
\begin{enumerate}[label=(\roman*)]
    \item There exists a right inverse $X_-^\dagger$ such that $X_+ X_-^\dagger$ is stable and 
    $$
    \begin{bmatrix} 
    W_- \\ Z_- 
    \end{bmatrix} X_-^\dagger = 0.
    $$ \label{cond1}
    \item There exist right inverses $X_-^\dagger$ and $W_-^\dagger$ such that $X_+ X_-^\dagger$ is stable, $W_- X_-^\dagger = 0$,
    $$
    \begin{bmatrix} 
    X_- \\ U_- 
    \end{bmatrix} W_-^\dagger = 0,
    $$
    and the unique solution $P$ to 
    \begin{equation}
    \label{Lyapunoveqdata}
       (X_-^\dagger)^\top \left( X_+^\top P X_+ - X_-^\top P X_- + Z_-^\top Z_- \right) X_-^\dagger = 0
    \end{equation}
    has the property that 
    \begin{equation}
    \label{traceeq}
        \trace \left( (X_+W_-^\dagger)^\top P X_+W_-^\dagger \right) < \gamma.
    \end{equation}
    \label{cond2}
\end{enumerate}
Moreover, $K$ is an $\mathcal{H}_2$ suboptimal controller for all $(A,B,E) \in \Sigma_{\ids}$ if and only if $K = U_- X_-^\dagger$, where $X_-^\dagger$ satisfies the conditions of \ref{cond1} or \ref{cond2}.
\end{thm}

\begin{rem}
The interpretation of Theorem \ref{thmH2} is as follows. Note that both condition \ref{cond1} and \ref{cond2} require the existence of $X_-^\dagger$ such that $X_+ X_-^\dagger$ is stable and $W_- X_-^\dagger = 0$. These conditions are necessary for the existence of a stabilizing controller by Lemma \ref{lemstabE}. In condition \ref{cond1} it is further required that $X_-^\dagger$ satisfies $Z_- X_-^\dagger = 0$, which means that the output of all systems in $\Sigma_{\ids}$ can be made identically equal to zero (hence the $\mathcal{H}_2$ norm is zero). In condition \ref{cond2}, the properties of $W_-^\dagger$ imply that $E_s = X_+ W_-^\dagger$ can be uniquely identified from the data. Similar to the suboptimal LQR problem, it is generally not required that $A_s$ and $B_s$ can be uniquely identified from the data.
\end{rem}

\begin{pf}
We first prove the `if' parts of both statements. Suppose that condition \ref{cond1} is satisfied and let $K := U_-X_-^\dagger$. By Lemma \ref{lemstabE}, $A+BK$ is stable for all $(A,B,E) \in \Sigma_{\ids}$. As $Z_- X_-^\dagger = 0$ we have $C + DU_-X_-^\dagger = C+DK = 0$. This means that the $\mathcal{H}_2$ norm of \eqref{closedloop} is zero for all $(A,B,E) \in \Sigma_{\ids}$. We conclude that the data are informative for $\mathcal{H}_2$ suboptimal control and $K$ is an $\mathcal{H}_2$ suboptimal controller. 

Next suppose that condition \ref{cond2} is satisfied, and let $K := U_- X_-^\dagger$ where $X_-^\dagger$ satisfies the conditions of \ref{cond2}. Clearly, $A+BK = X_+X_-^\dagger$ is stable for all $(A,B,E) \in \Sigma_{\ids}$. By the properties of $W_-^\dagger$, $(A,B,E) \in \Sigma_{\ids}$ implies $E = E_s$. In view of \eqref{Lyapunoveqdata} and \eqref{traceeq} we see that for any $(A,B,E_s) \in \Sigma_{\ids}$ the unique solution $P$ to \eqref{Lyapunoveq} satisfies $\trace (E_s^\top P E_s) < \gamma$. Therefore, the data are informative for $\mathcal{H}_2$ suboptimal control and $K$ is $\mathcal{H}_2$ suboptimal. 

Subsequently, we prove the `only if' parts of both statements. Suppose that the data are informative for $\mathcal{H}_2$ suboptimal control and let $K$ be an $\mathcal{H}_2$ suboptimal controller for all $(A,B,E) \in \Sigma_{\ids}$. By Lemma \ref{lemstabE}, there exists a right inverse $X_-^\dagger$ such that $X_+X_-^\dagger$ is stable and $W_-X_-^\dagger = 0$. Also, the feedback $K$ is of the form $K = U_- X_-^\dagger$. The solution $P$ to \eqref{Lyapunoveqdata} exists and is unique by stability of $X_+X_-^\dagger$. The matrix $P$ satisfies $\trace(E^\top P E) < \gamma$ for all $(A,B,E) \in \Sigma_{\ids}$. Therefore, we have 
\begin{equation}
\label{E0}
    \trace\left((E+\alpha E_0)^\top P (E+\alpha E_0)\right) < \gamma
\end{equation}
for all $(A,B,E) \in \Sigma_{\ids}$, $(A_0,B_0,E_0) \in \Sigma^0_{\ids}$ and $\alpha \in \mathbb{R}$. We divide both sides of \eqref{E0} by $\alpha^2$ and take the limit as $\alpha \to \infty$. Then, by continuity of the trace we obtain $\trace(E_0^\top P E_0) = 0$, which yields $P E_0 = 0$ for all $(A_0,B_0,E_0) \in \Sigma^0_{\ids}$. We claim that this implies that either $P = 0$ or $E_0 = 0$ for all $(A_0,B_0,E_0) \in \Sigma^0_{\ids}$. Suppose that this claim is not true. Then $P \neq 0$ and there exists a triple $(A_0,B_0,E_0) \in \Sigma^0_{\ids}$ such that $E_0 \neq 0$. Note that $(F A_0,F B_0,F E_0) \in \Sigma^0_{\ids}$ for any $F \in \mathbb{R}^{n \times n}$. Clearly, there exists an $F$ such that $P F E_0 \neq 0$. This is a contradiction, which proves our claim. Now, in the case that $P = 0$ we obtain $Z_- X_-^\dagger$ and condition \ref{cond1} is satisfied. In the case that $E_0=0$ for all $(A_0,B_0,E_0) \in \Sigma^0_{\ids}$, there exists a right inverse $W_-^\dagger$ such that $X_- W_-^\dagger = 0$ and $U_- W_-^\dagger = 0$. This means that $(A,B,E) \in \Sigma_{\ids}$ implies $E = E_s = X_+ W_-^\dagger$. Hence \eqref{traceeq}, and therefore \ref{cond2}, holds. In both cases, the controller $K$ is of the form $K = U_- X_-^\dagger$, where $X_-^\dagger$ satisfies either \ref{cond1} or \ref{cond2}. \hfill $\square$
\end{pf}
Similar to Corollary \ref{corsuboptimalLQR} we can reformulate Theorem \ref{thmH2} in terms of linear matrix inequalities using Proposition \ref{propsuboptimalH2}. 

\begin{cor}
Let $\gamma > 0$. The data $(U_-,W_-,X)$ are informative for $\mathcal{H}_2$ suboptimal control if and only if at least one of the following two conditions is satisfied:
\begin{enumerate}[label=(\roman*)]
    \item There exists a $\Theta \in \mathbb{R}^{T \times n}$ such that $X_- \Theta = (X_-\Theta)^\top$,
    $$
    \begin{bmatrix} 
    W_- \\ Z_- 
    \end{bmatrix} \Theta = 0 \text{ and } \begin{bmatrix}
    X_- \Theta & \Theta^\top X_+^\top \\ X_+ \Theta & X_-\Theta 
    \end{bmatrix} > 0.
    $$ \label{corcond1}
    \item There exists a right inverse $W_-^\dagger$, a $Y = Y^\top \in \mathbb{R}^{n \times n}$ and $\Theta \in \mathbb{R}^{T \times n}$ such that $X_-\Theta$ is symmetric, the matrices $W_- \Theta$, $X_- W_-^\dagger$ and $U_- W_-^\dagger$ are zero, and
    \begin{align*}
        \begin{bmatrix}
        X_-\Theta & \Theta^\top X_+^\top & \Theta^\top Z_-^\top \\ X_+ \Theta & X_-\Theta & 0 \\ Z_- \Theta & 0 & I
        \end{bmatrix} &> 0 \\
        \begin{bmatrix}
        Y & (W_-^\dagger)^\top X_+^\top \\ X_+ W_-^\dagger & X_-\Theta
        \end{bmatrix} &> 0 \\
        \trace(Y) &< \gamma.
    \end{align*} \label{corcond2}
    \end{enumerate}
    Moreover, $K$ is an $\mathcal{H}_2$ suboptimal controller for all $(A,B,E) \in \Sigma_{\ids}$ if and only if $K = U_- \Theta (X_-\Theta)^{-1}$, where $\Theta$ satisfies the conditions of \ref{corcond1} or \ref{corcond2}.
\end{cor}

\section{Illustrative example}
\label{sectionexample}

\begin{figure}[h!]
\centering
\includegraphics[width=0.33\textwidth]{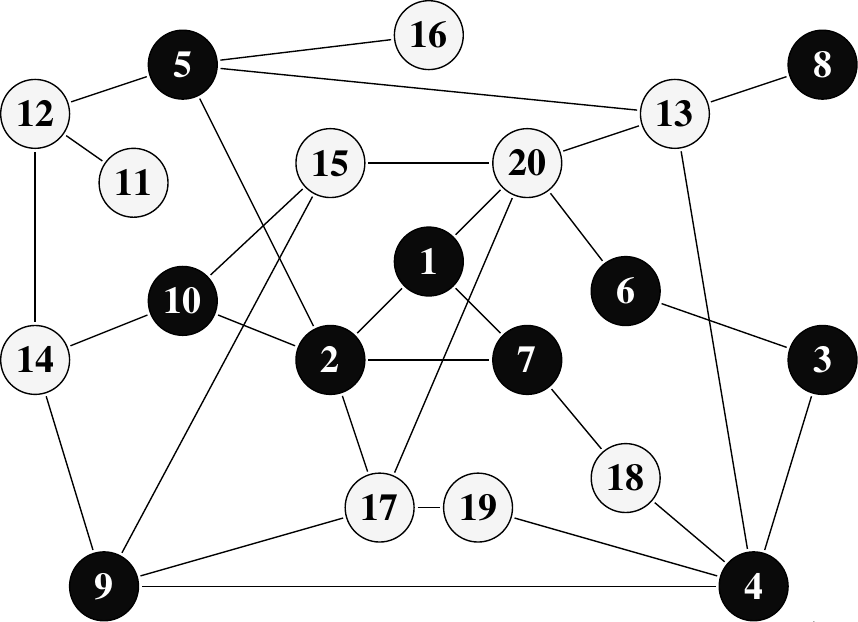}
\caption{Graph $G$ with leader vertices colored black.}
\label{fig:graph}
\end{figure}

We study steered consensus dynamics of the form
\begin{equation}
\label{consensus}
    \bmx(t+1) = \left(I - 0.15 L\right) \bmx(t) + B \bmu(t),
\end{equation}
where $\bmx \in \mathbb{R}^{20}$, $\bmu \in \mathbb{R}^{10}$, $L$ is the Laplacian matrix of the graph $G$ in Figure \ref{fig:graph}, and $B = \begin{bmatrix}
I & 0
\end{bmatrix}^\top$, meaning that inputs are applied to the first 10 nodes. The goal of this example is to apply the theory from \S \ref{sectionLQR} to construct suboptimal controllers for \eqref{consensus} using data. We choose the weight matrices as $Q = I$ and $R = I$, and define $x_0 \in \mathbb{R}^{20}$ entry-wise as $(x_0)_i = i$.

We start with a time horizon of $T = 20$ and collect data $(U_-,X)$ where the entries of $U_-$ and the initial state of the experiment $x(0)$ are drawn uniformly at random from $(0,1)$. Given these data, we attempt to solve a semidefinite program (SDP) where the objective is to minimize $\gamma$ subject to the constraints \eqref{1}, \eqref{2} and \eqref{3}. We use Yalmip, with Mosek as a solver. Next, we collect one additional sample of the input and state, and we solve the SDP again for the augmented data set. We continue this process up to a time horizon of $T = 30$. 

\begin{figure}[h!]
\centering
\includegraphics[width=0.48\textwidth,height=5cm]{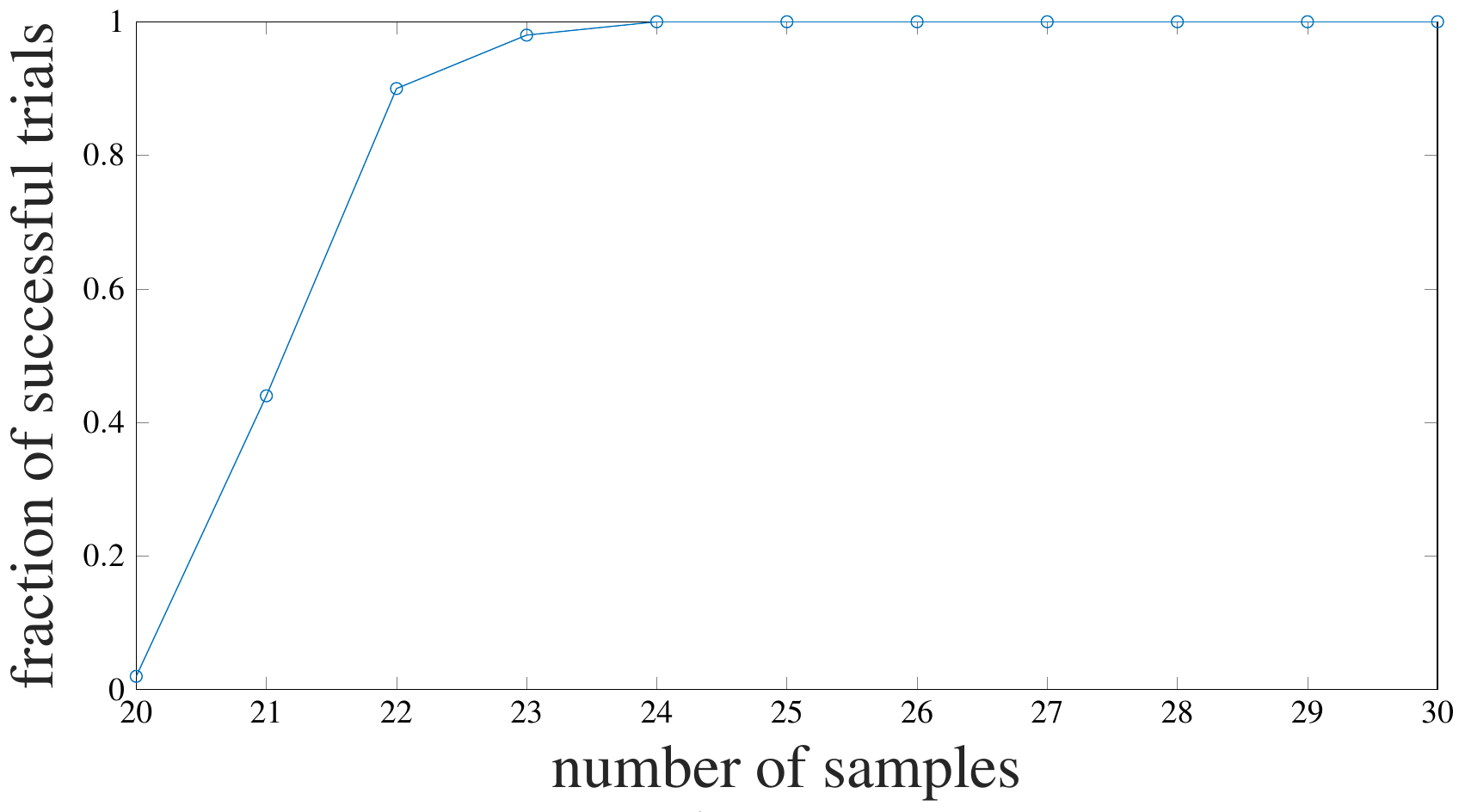}
\caption{Fraction of successful trials as a function of $T$.}
\label{fig:fraction}
\end{figure}

\begin{figure}[h!]
\centering
\includegraphics[width=0.48\textwidth,height=5cm]{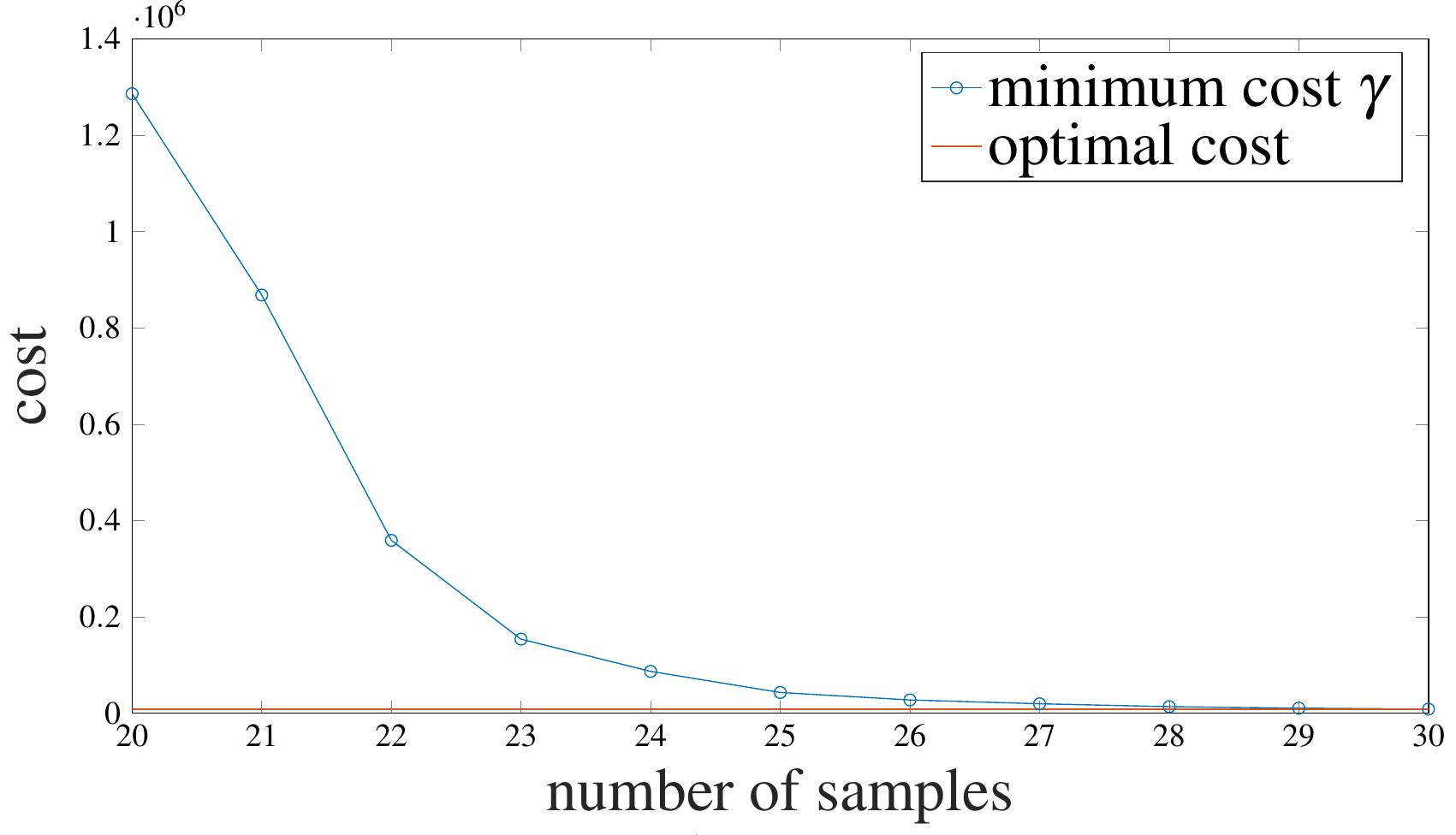}
\caption{Average minimum cost as a function of $T$.}
\label{fig:cost}
\end{figure}

We repeat this entire experiment for 100 trials and display the results in Figures \ref{fig:fraction} and \ref{fig:cost}. Figure \ref{fig:fraction} depicts the fraction of successful trials in which the constraints \eqref{1}, \eqref{2} and \eqref{3} were feasible and a stabilizing controller was found. Note that a stabilizing controller was only found in 2 out of the 100 trials for $T=20$. This fraction rapidly increases to $0.88$ for $T = 22$, while $100\%$ of the trials were successful for $T \geq 24$. Figure \ref{fig:cost} displays the minimum cost $\gamma$ of the controller, averaged over all successful trials. The cost is very large for small sample size $(T = 20)$ but decreases rapidly as the number of samples increases. Figure \ref{fig:cost} therefore highlights an interesting trade-off between the sample size and the cost. Note that for $T = 30$, $\gamma$ coincides with the optimal cost obtained from the (model-based) solution to the Riccati equation. This is as expected since $30 = n+m$ is the minimum number of samples from which the state and input matrices can be uniquely identified.

\section{Conclusions}
\label{sectionconclusions}
In this paper we have studied the data-driven suboptimal LQR and $\mathcal{H}_2$ problems. For both problems, we have presented conditions under which a given data set contains sufficient information for control design. We have also given a parameterization of all suboptimal controllers in terms of data-driven linear matrix inequalities. Finally, we have illustrated these results by numerical simulations, which reveal a trade-off between the number of collected data samples and the achieved controller performance. 

\bibliography{ifacconf}             

\end{document}